\documentclass[12pt,oneside,a4paper]{amsart}

\usepackage{amscd}

\usepackage{amscd}
\usepackage{amssymb}
\usepackage{a4wide}
\usepackage{amstext}
\usepackage{amsthm}
\usepackage{mathrsfs}

\usepackage[final]{hyperref}
\hypersetup{unicode= false, colorlinks=true, linkcolor=blue,
anchorcolor=blue, citecolor=green, filecolor=red, menucolor=blue,
pagecolor=blue, urlcolor=blue} 
\numberwithin{equation}{section}

\newcommand{\CC}{\mathbb {C}}

\newcommand{\RR}{\mathbb{R}}

\DeclareMathOperator{\ospan}{\overline{Span}}
\DeclareMathOperator{\sig}{sign}

\renewcommand{\phi}{\varphi}
\newcommand{\pw}{\mathcal{P}W_\pi}

\newcommand{\co}{\mathbb{C}}
\newcommand{\na}{\mathbb{N}}
\newcommand{\zl}{\mathbb{Z}}

\newcommand{\rl}{\mathbb{R}}

\newcommand{\he}{\mathcal{H}(E)}

\newcommand{\la}{\lambda}
\newcommand{\BC}{\mathbb{C}}
\newcommand{\BR}{\mathbb{R} }

\renewcommand{\phi}{\varphi}
\newcommand{\deab}{\varkappa}

\newcommand{\cD}{\mathcal{D}}

\newcommand{\cDA}{\cD(A)}
\newcommand{\cDL}{\cD(L)}

\newcommand{\Ker}{\operatorname{Ker}}

\newtheorem{Thm}{Theorem}[section]
\newtheorem{theorem}[Thm]{Theorem}

\newtheorem{corollary}[Thm]{Corollary}

\newtheorem{example}[Thm]{Example}

\begin{document}
\sloppy
\title[Recent developments in spectral synthesis]
{Recent developments in spectral synthesis for exponential
systems and for non-self-adjoint operators}

\author{Anton Baranov, Yurii Belov, Alexander Borichev, Dmitry Yakubovich}
\address{Anton Baranov,
\newline Department of Mathematics and Mechanics, St. Petersburg State University, St. Petersburg, Russia,
\newline {\tt anton.d.baranov@gmail.com}
\smallskip
\newline \phantom{x}\,\, Yurii Belov,
\newline Chebyshev Laboratory, St. Petersburg State University, St. Petersburg, Russia,
\newline {\tt j\_b\_juri\_belov@mail.ru}
\smallskip
\newline \phantom{x}\,\, Alexander Borichev,
\newline Aix-Marseille Universit\'e, Marseille, France
\newline {\tt borichev@cmi.univ-mrs.fr}
\smallskip
\newline \phantom{x}\,\, Dmitry Yakubovich,
\newline Departamento de Matem\'{a}ticas, Universidad Autonoma de Madrid, Spain
\newline {\tt dmitry.yakubovich@uam.es}
}
\thanks{A. Baranov and Yu. Belov were supported by the Chebyshev Laboratory
(St. Petersburg State University) under the Russian Federation Government grant 11.G34.31.0026.
A. Baranov was partially supported by RFBR grant 11-01-00584-a.
The research of A. Borichev was partially supported by
the ANR FRAB.
D. Yakubovich has been supported by the
Project MTM2011-28149-C02-01
and by the ICMAT Severo Ochoa project SEV-2011-0087
of the Ministry of Economy and Competition of Spain.}

\begin{abstract}
We survey recent results concerning the hereditary completeness
of some special systems of functions and the spectral
synthesis problem for a related class of linear operators.
We present a solution of the spectral synthesis problem
for systems of exponentials in $L^2(-\pi, \pi)$.
Analogous results are obtained for the systems of reproducing
kernels in the de Branges spaces of entire functions.
We also apply these results (via a functional model)
to the spectral theory of rank one perturbations of compact self-adjoint
operators.
\end{abstract}

\dedicatory{Dedicated to Nikolai K. Nikolski on the occasion of his 70-th birthday}

\maketitle

\section{Introduction}


Spectral synthesis for linear operators and hereditary completeness of systems of vectors
in a Hilbert space are among numerous mathematical interests of Nikolai Nikolski.
His influence on the subject is enormous.
The authors are deeply grateful to Nikolai Nikolski,
who introduced them to the field of the spectral function theory
and posed many of the problems studied in the paper.

The aim of the present
paper is to give a short and accessible overview of our recent results on two closely
related topics. The first one concerns the hereditary completeness for
systems of exponentials in $L^2(-\pi, \pi)$ or, more generally,
for systems of reproducing kernels in the de Branges spaces of
entire functions (studied by the first three authors in
\cite{bbb}). The second topic is the spectral theory and, in
particular, the spectral synthesis problem for rank one
perturbations of compact self-adjoint operators. This part is based
on the preprints \cite{bar-yak, bar-yak1} of the first and the
fourth authors.


\subsection{Hereditary completeness}
A system of vectors  $\{x_n\}_{n\in N}$ in a separable Hilbert space
$H$ is said to be {\it exact} if it is both {\it complete}
(i.e., $\ospan \{x_n\} = H$)
and {\it minimal}
(i.e., $\ospan \{x_n\}_{n\neq n_0} \neq H$ for any $n_0$).
Given an exact system we consider its (unique) {\it biorthogonal}
system $\{\tilde x_n\}_{n\in N}$ which satisfies
$ (x_m, \tilde x_n) = \delta_{mn}$.
We associate to every element $x\in H$ its formal Fourier series
\begin{equation}
\label{rf}
x  \sim \sum_{n\in N}(x,\tilde x_n)x_n.
\end{equation}

In what follows, we will also need the notion of a \textit{linear summation method} for
a vector series
$\sum_{n\ge 1} v_n$,
where $v_n$ are vectors in $H$. Let $R=(R_{jk})_{j,k=1}^\infty$
be a doubly infinite matrix, whose entries satisfy the consistency properties
$$
\sup_{j\ge 1}\sum_{k=1}^\infty|R_{jk}|<\infty,\qquad \lim_{j\to\infty}R_{jk}=0,\quad k\ge 1,\qquad
\lim_{j\to\infty}\sum_{k=1}^\infty R_{jk}=1.
$$
The series $\sum_{n\ge 1} v_n$ is summable to a vector $v$ by the method $R$ if the limit
$$
\lim_{j\to\infty}\sum_{k=1}^\infty R_{jk}\sum_{n=1}^k v_n
$$
exists and equals to $v$.


The system $\{x_n\}_{n\in N}$ is said to be
{\it hereditarily complete} if, for any $x\in H$,
$$
x\in \ospan \{(x,\tilde x_n) x_n\},
$$
which may be understood as a possibility of reconstructing
(in a very weak sense) the vector $x$ from its Fourier series.
Clearly, if
the series (\ref{rf}) admits a
linear summation method,
then the system $\{x_n\}$ is hereditarily complete.

An equivalent definition of hereditary completeness is that
for any partition $N = N_1 \cup N_2$, $N_1 \cap N_2 =\emptyset$,
the system
$$
\{x_n\}_{n\in N_1} \cup \{\tilde x_n\}_{n\in N_2}
$$
is complete in $H$. In particular, this requires that
the biorthogonal system $\{\tilde x_n\}$ be also complete (which is not always true).
If both $\{x_n\}$ and $\{\tilde x_n\}$ are complete,
but the system is not hereditary complete, we say
that $\{x_n\}$ is {\it nonhereditarily complete}.

A general approach to constructing nonhereditarily
complete systems was developed by Dovbysh, Nikolski and Sudakov
\cite{dn, dns}, where it was shown, in particular, that
for any condition of closeness to an orthogonal basis
(weaker than the quadratical closeness)
there are uniformly minimal nonhereditarily complete systems
that meet this condition.
We recall that a complete minimal system that is quadratically
close to an orthogonal basis is always a Riesz basis, by a classical
theorem of Bari.
%
%
%


\subsection{Spectral synthesis of linear operators}
Let $T$ be a bounded linear operator in a Hilbert space $H$
which has a complete set of eigenvectors (or, more generally, a complete set of root vectors) $\{x_n\}$.
The operator $T$ is said to admit the {\it spectral synthesis} if
for any $T$-invariant closed linear subspace $M$ of $H$
the restriction of $T$ onto $M$ has a complete set of root vectors.
Equivalently, $M = \ospan \{x_n: x_n \in M\}$, that is, $M$
is generated by the root vectors it contains.

The spectral synthesis problem in this form was posed by Wermer \cite{wer}
who showed that any unitary or self-adjoint operator with a complete set of
eigenvectors admits the spectral synthesis as well as any compact normal
operator. The first example of a compact operator
with a complete set of
eigenvectors (and with a trivial kernel) which does not admit the spectral
synthesis goes back to Hamburger \cite{hamb}, who constructed
a compact operator with a complete set of eigenvectors, whose restriction
to an invariant subspace is a nonzero Volterra operator.

Further results about the spectral synthesis
of compact operators were obtained by
Nikolski \cite{nik69} and Markus \cite{markus}. Let us mention
the following theorem of Nikolski
\cite{nik69}, which shows that any Volterra operator
can be realized as a part of compact operator with a complete set of eigenvectors:
{\it For any Volterra operator $T$ in a Hilbert space $H$ there exists a
Hilbert space $H'$ and a compact operator $T': H\oplus H' \to H\oplus H'$
such that $T'H\subset H$, $T'|_H = T$, and the eigenfunctions of $T'$
corresponding to nonzero eigenvalues are complete in $H\oplus H'$. }
In this situation one can control the decay of the singular values of $T'$
in terms of the corresponding information on $T$, see \cite{Solom-B}.

On the other hand, it was shown in \cite{markus} that
if $T$ is a compact operator
with a complete set of eigenvectors $\{x_n\}$, then
$T$ admits the spectral synthesis if and only if
$\{x_n\}$ is hereditarily complete.
For an extensive review of the subject see the surveys \cite{niko1}, \cite[Chapter 4]{hrnik}, and \cite{niksel}.


\subsection{Hereditary completeness for exponential systems}
It is natural to study the problem of hereditary completeness
for special systems in functional spaces, e.g. those which appear
as families of eigenvectors and root vectors of a certain
operator. Exponential systems form an important class in this respect.
Let $\Lambda = \{\lambda_n\} \subset \co$ and let
$e_\lambda (t) = \exp(i\lambda t)$. We consider
the exponential system $\{e_\lambda\}_{\lambda\in\Lambda}$ in $L^2(-a, a)$, $a>0$.
By the result of Young \cite{young},
for any exact system of exponentials $\{e_\lambda\}_{\lambda\in\Lambda}$
its biorthogonal system
$\{\tilde e_\lambda\}_{\lambda\in\Lambda}$ is complete.

Another important class of systems are families of the reproducing kernels
in the so-called model subspaces $K_\Theta = H^2\ominus\Theta H^2$
of the Hardy space $H^2$ (either in the disk or in the half plane), where
$\Theta$ is an inner function. The spaces $K_\Theta$ are known
especially for their role in the Nagy--Foia\c{s} functional model. They
have numerous applications in operator-related function theory
(we refer to Nikolski's monographs \cite{nik0, nk12}).

The following two questions were posed by Nikolski:
%
%
\medskip

{\it Does there exist an exact system of reproducing kernels in $K_\Theta$\textup,
such that its biorthogonal system is incomplete\textup?}
\smallskip

{\it Does there exist an exact system of reproducing kernels in $K_\Theta$
with the complete biorthogonal system\textup, which is nonhereditarily complete\textup? }
\medskip

Both questions have positive answers.
Examples of exact systems of reproducing kernels in $K_\Theta$
with incomplete biorthogonal systems were constructed
in \cite{bb} (see Theorem \ref{biort}).

There is a unitary operator (constructed from the Fourier
transform), that sends any exponential system in
$L^2(-a,a)$ to
a system of reproducing kernels
in the model space $K_{e^{2iaz}}$ in the upper half-plane.
Therefore the hereditary completeness of exponential systems is a special case
of the second question. One of the main results of \cite{bbb}
is that there exist nonhereditarily complete systems of exponentials.

\begin{theorem}
\label{maincount0}
There exists a system of exponentials $\{e^{i\lambda t}\}_{\lambda\in \Lambda}$\textup,
$\Lambda \subset \rl$\textup,
which is complete and minimal in $L^2(-\pi, \pi)$\textup, but
is not hereditarily complete.
\end{theorem}

Thus, in general there exists no linear summation method
for nonharmonic Fourier series.

\begin{corollary}
\label{maincount00}
There exist a system of exponentials $\{e^{i\lambda t}\}_{\lambda\in \Lambda}$\textup,
$\Lambda \subset \rl$\textup, which is complete and minimal in $L^2(-\pi, \pi)$\textup,
such that the corresponding Fourier series $\sum_{\lambda\in \Lambda}
(f,\tilde e_\lambda ) e_\lambda$ admit no linear summation method.
\end{corollary}

However, surprisingly,
the exponential systems are hereditarily complete up to
a one-dimensional defect (see Theorem \ref{main1}
and  the discussion in the next section).
Analogous results are obtained for systems of reproducing kernels
in the de Branges spaces of entire functions (which correspond to
model spaces generated by inner functions $\Theta$
meromorphic in $\mathbb C$), see Section 3.


\subsection{One-dimensional perturbations of compact self-adjoint
operators}
Systems of reproducing kernels in Hilbert spaces of analytic
functions may serve as models for abstract systems of vectors
in a general Hilbert space. As was pointed out by Nikolski,
the model subspaces form a promising setting in this respect
due to their rich structure.

One may ask, whether
the exact systems of exponentials or
reproducing kernels of de Branges spaces may be interpreted
as eigenfunctions of bounded linear operators
from some special class. Here the
rank one perturbations of bounded self-adjoint operators
come into play.

Let $A^\circ$ be a compact self-adjoint operator on a Hilbert space $H$
with simple spectrum and with a trivial kernel.
Let $a^\circ, b^\circ \in H$ and let
$L^\circ = L(A^\circ, a^\circ, b^\circ)$
be its one-dimensional perturbation\footnote{We use the notation $L^\circ$, $A^\circ$ etc.
to distinguish the compact operators from the case of {\it singular }
rank one perturbations of {\it unbounded} self-adjoint operators,
which will be respectively denoted by
$L$ and $A$.}, that is,
$$
L^\circ f = A^\circ f + (f, b^\circ) a^\circ, \qquad f\in H.
$$

%
Though
this is one of the simplest and most natural classes of compact linear operators, it
seems that our knowledge of the
spectral structure of these operators is very incomplete.

One of the basic questions addressed in  \cite{bar-yak} was
to understand when the perturbed operator $L^\circ$
has a complete set of eigenfunctions. Notice that the classical theory
covers only the case of so-called weak perturbations, i.e., perturbations
of the form $L^\circ = A^\circ(I+T^\circ)$ where $T^\circ$ is a compact operator
(theorems of Keldysh and Matsaev), as well as some other special classes
of operators (e.g., dissipative ones), see \cite{Mats61} and
\cite[Chapter V]{Gohb_Krein}.

In our study, we use a functional model
for rank-one perturbations of self-adjoint operators, which is close to a model by
Kapustin given in \cite{kap}.
See Section 4 below for more details. This model allows us to obtain
new positive results on completeness and spectral synthesis
for rank one perturbations of compact self-adjoint operators
as well as a number of rather unexpected counterexamples.
These counterexamples show how rich and complicated the spectral
structure becomes even for rank one perturbations as soon as we leave the
area of dissipative operators.

Here we mention one of these examples
(for the precise statement see Theorem \ref{synthesis}):

\begin{theorem}
\label{synthesis-1}
Under some mild separation conditions on the spectrum
of $A^\circ$ there exists a rank one perturbation $L^\circ$
of $A^\circ$ which does not admit the spectral synthesis.
\end{theorem}
\medskip

The paper is organized as follows. In Section 2
we discuss the solution of the spectral synthesis problem for exponential
systems, while in Section 3 some extensions to the general
de Branges spaces are considered.
The functional model for rank one perturbations and its corollaries
are discussed in Sections 4 and 5. We conclude the article
with a list of open problems.


\section{Hereditary completeness for exponential systems}

As usual, applying the Fourier transform $\mathcal{F}$
one reduces the problem for exponential systems
in $L^2(-\pi, \pi)$
to the same problem for systems of reproducing kernels in
the Paley--Wiener space
$\pw = \mathcal{F} L^2(-\pi,\pi)$. Recall
that the reproducing kernel of $\pw$ corresponding to a point $\lambda\in\co$
is of the form
$$
K_\lambda(z) = \frac{\sin \pi (z-\overline \lambda)}{\pi(z-\overline
\lambda)}, \qquad f(\lambda) = (f,K_\lambda)_{\pw}.
$$

In \cite{bbb}, we solve the problem of hereditary completeness
for exponential systems. Namely, we show that  the hereditary completeness
holds up to a possible one-dimensional defect.

Let $\Lambda\subset \CC$ be such that the system of reproducing kernels
$\{K_\lambda\}_{\lambda\in \Lambda}$ is exact
in the Paley--Wiener space $\pw$.
Then the biorthogonal system $\{g_\lambda\}$
is given by
$$
g_\lambda(z) = \frac{G(z)}{G'(\lambda)(z-\lambda)},
$$
where $G$ is the so-called {\it generating function} of the set
$\Lambda$: the entire function $G$ has only simple zeros at the points of $\Lambda$ and
$g_\lambda\in \pw$, $\lambda\in\Lambda$.
By the above-mentioned result of Young,
$\{g_\lambda\}$  is also an exact system.
It is well known that $G$ is 
of exponential
type $\pi$.

The following theorem gives a complete answer to the
hereditary completeness of reproducing kernels
(equivalently exponential systems).

\begin{theorem}[\cite{bbb}, Theorems 1.1, 1.3]
\label{main1}
1. If $\{K_\lambda\}_{\lambda\in \Lambda}$ is exact
in the Paley--Wiener space $\pw$\textup,
then for any partition $\Lambda = \Lambda_1 \cup \Lambda_2$\textup, the
orthogonal complement in $\pw$ to the system
\begin{equation}
\label{syst}
\{g_\lambda\}_{\lambda\in \Lambda_1} \cup
\{K_\lambda\}_{\lambda\in \Lambda_2}
\end{equation}
is at most one-dimensional.

2. There exists a system $\{K_\lambda\}_{\lambda\in \Lambda}$\textup,
$\Lambda \subset \rl$\textup, which is exact in $\pw$\textup, but
is not hereditarily complete \textup(thus\textup, there exists a partition
$\Lambda = \Lambda_1 \cup \Lambda_2$\textup, such that
the orthogonal complement in $\pw$ to the system \eqref{syst}
is one-dimensional\textup).
\end{theorem}

Thus, hereditary completeness may fail even for exponential systems
(reproducing kernels of the Paley--Wiener space), which answers
the question of Nikolski.  Further counterexamples will be discussed in
the next section.

The set of the exceptional partitions
$\Lambda = \Lambda_1\cup\Lambda_2$ (for which the orthogonal
complement to \eqref{syst} is nontrivial) is in a sense very
small. This nontrivial defect cannot appear
unless the sequence $\Lambda_1$ has zero density. Given a sequence $\Lambda$, set
$$
D_+(\Lambda) = \limsup_{r\to\infty} \frac{n_r(\Lambda)}{2r},
$$
where $n_r(\Lambda)$ is the usual counting function of the
sequence $\Lambda$, $n_r(\Lambda) =\#\{ {\lambda\in\Lambda, |\lambda|\leq r}\}$.

\begin{theorem}[\cite{bbb}, Theorem 1.2]
\label{main2}
Let $\Lambda \subset \mathbb{C}$\textup, let the system $\{K_\lambda\}_{\lambda\in \Lambda}$
be exact in $\pw$\textup, and let the partition $\Lambda = \Lambda_1 \cup \Lambda_2$
satisfy $D_+(\Lambda_1) >0$.
Then the system \eqref{syst} is complete in $\pw$.
\end{theorem}

Theorem \ref{main2} shows that there is a strong asymmetry
between the systems of reproducing kernels and their biorthogonal.
The completeness of (\ref{syst}) may fail only when we take
a sparse (but infinite!) subsequence $\Lambda_1$.


\section{Systems of reproducing kernels in the de Branges spaces}

\subsection{Preliminaries on de Branges spaces}
Let $E$ be an entire function in the Hermite--Biehler class, that is
$E$ has no zeros in $\mathbb{C}_+ \cup\rl$, and
$$
|E(z)| > |E^*(z)|, \qquad z\in {\mathbb{C}_+},
$$
where $E^* (z) = \overline {E (\overline z)}$.
With any such function we associate the {\it de Branges space}
$\mathcal{H} (E) $ which consists of all entire functions
$F$ such that $F/E$ and $F^*/E$ restricted to $\mathbb{C_+}$ belong
to the Hardy space $H^2=H^2(\mathbb{C_+})$.
The inner product in $\he$ is given by
$$
( F,G)_E = \int_\rl \frac{F(t)\overline{G(t)}}{|E(t)|^2} \,dt.
$$
The reproducing kernel of the de Branges space ${\mathcal H} (E)$
corresponding to the point $w\in \mathbb{C}$ is given by
$$
K_w(z)=\frac{\overline{E(w)} E(z) - \overline{E^*(w)} E^*(z)}
{2\pi i(\overline w-z)} =
\frac{\overline{A_E(w)} B_E(z) -\overline{B_E(w)}A_E(z)}{\pi(z-\overline w)},
$$
where we use the standard decomposition $E=A_E - iB_E$, $A_E = \frac{E+E^*}{2}$,
$B_E=\frac{E^*-E}{2i}$.

The Hilbert spaces of entire functions $\he$ were introduced by
L. de Branges \cite{br} in relation to the inverse
spectral problems for differential operators.
These spaces are also of a great interest from the function theory
point of view. The Paley--Wiener space $\mathcal{P}W_a$
is the de Branges space corresponding to
$E(z) = \exp(-iaz)$. Also, de Branges spaces
are canonically isomorphic to
model spaces generated by meromorphic inner functions:
if $\Theta$ has a meromorphic continuation
to the whole plane, then $\Theta = E^*/E$ for a function $E$
in the Hermite--Biehler class and the mapping $f\mapsto Ef$
is a unitary operator from $K_\Theta$ onto $\mathcal{H}(E)$,
which maps reproducing kernels onto reproducing kernels.

An important characteristics of the de Branges space $\he$
is its phase function, that is, an increasing $C^\infty$-function
$\phi$ such that $E(t) \exp{i \phi(t)} \in \rl$, $t\in\rl$
(thus, essentially, $\phi = - \arg E$ on $\rl$).
Clearly, for $\he=\mathcal{P}W_a$, $\phi(t) = at$.
If $\phi' \in L^\infty(\rl)$ (in which case we say that $\phi$
has sublinear growth), the space $\he$ shares some properties
with the Paley--Wiener spaces.

A crucial property of the de Branges spaces is the existence of orthogonal bases
of reproducing kernels corresponding to real points \cite{br}.
For $\alpha\in [0, \pi)$
we consider the set of points $t_n\in\mathbb{R}$
such that
\begin{equation}
\label{basis}
\varphi(t_n)=\alpha+\pi n,\qquad n\in\mathbb{Z}.
\end{equation}
Thus, $\{t_n\}$ is the zero set of the function
$e^{i\alpha} E - e^{-i\alpha}E^*$.

If the points $t_n$ are defined by (\ref{basis}),
then the system of reproducing kernels $\{K_{t_n}\}$
is an orthogonal basis for $\he$ for each $\alpha\in [0, \pi)$
except, may be, one ($\alpha$ is an exceptional value if and only if
$e^{i\alpha} E - e^{-i\alpha}E^* \in \he$).
In particular, if $A_E \notin \he$, then $\big\{\frac{A_E}{z-t_n}\big\}$
is the orthogonal basis for $\he$ corresponding to $\alpha = \pi\slash2$.

It should be mentioned that the index set for the sequence $\{t_n\}$ is either
$\mathbb Z$ or $\mathbb Z_\pm$; in the latter case we may need to add $k\pi$ to $\varphi$, for some $k\in\mathbb Z$.
%
%

The norm of $K_t$, $t\in \rl$, is given by
$\|K_t\|_E^2 = \pi^{-1} |E(t)|^2 \phi'(t)$.
Thus, if we put $\mu = \sum_n \|K_{t_n}\|^{-2} \delta_{t_n}$,
then the embedding $\he \to L^2(\mu)$ is a unitary operator.
One should think of the sequence $\{t_n\}$
as of a spectral characteristics of the space $\he$.


\subsection{Hereditary completeness for reproducing kernels of
the de Branges spaces}
As in the Paley--Wiener case we denote by
$G$ the generating function
of the set $\Lambda$, and
the system biorthogonal to the system
$\{K_\lambda\}_{\lambda\in \Lambda}$
is given by $\{g_\lambda\}_{\lambda\in \Lambda}$,
where $g_\lambda(z) = \frac{G(z)}{G'(\lambda)(z-\lambda)}$.

We say that $\phi$ is of {\it tempered growth} if $\phi'(t)=O(|t|^N)$, $|t|\to\infty$, for some $N$.

The completeness of the systems biorthogonal to exact systems
of reproducing kernels was studied in \cite{bb, fric}.
In particular, it was shown in \cite{fric}
that such biorthogonal systems are always
complete when $\phi'\in L^\infty(\rl)$.

The following result is obtained in \cite{bb}:

\begin{theorem}[\cite{bb}, Theorems 1.1, 1.2]
\label{biort}
1. If $e^{i\alpha} E - e^{-i\alpha}E^* \in \he$ for some $\alpha \in [0,\pi)$\textup,
then there exists an exact system of reproducing kernels\textup,
whose biorthogonal is incomplete.

2. If $\phi$ is of tempered growth\textup, then
any system $\{g_\lambda\}_{\lambda\in \Lambda} $
biorthogonal to an exact system of reproducing kernels
has at most a finite-dimensional defect in $\he$.
If\textup, moreover\textup, $e^{i\alpha} E - e^{-i\alpha}E^* \notin \he$ for any
$\alpha \in [0,\pi)$\textup, then $\{g_\lambda\}_{\lambda\in \Lambda} $
is complete in $\he$.
\end{theorem}

The method of the proof of Theorem \ref{main1}
extends to the case of the de Branges spaces with sublinear
or tempered growth of the phase.

\begin{theorem}[\cite{bbb}, Theorem 1.4, Theorem 5.3]
\label{main3}
Let $\he$ be a de Branges space such that
$\phi$ is of tempered growth.
If the system of reproducing kernels
$\{K_\lambda\}_{\lambda\in \Lambda}$ is exact in $\he$\textup,
then for any partition $\Lambda = \Lambda_1 \cup \Lambda_2$\textup,
the orthogonal complement in $\he$ to the system
\begin{equation}
\label{syst3}
\{g_\lambda\}_{\lambda\in \Lambda_1} \cup
\{K_\lambda\}_{\lambda\in \Lambda_2}
\end{equation}
is of finite dimension.

If\textup, moreover\textup, $\phi'\in L^\infty(\rl)$\textup, then
the orthogonal complement in $\he$ to the system \eqref{syst3}
is at most one-dimensional.
\end{theorem}

A crucial step in the proofs of Theorems \ref{main1} and \ref{main3} is the use
of expansions of functions
in $\pw$ or in $\he$ with respect to {\it two different}
orthogonal bases of reproducing kernels. At first glance
it may look like an artificial trick; however it should be noted that
the existence of two orthogonal bases of reproducing kernels
is a property which characterizes de Branges spaces
among all Hilbert spaces of entire functions (see \cite{BMS, BMS1}).
Therefore, we believe this method to be intrinsically related to the deep and complicated geometry of de Branges spaces.

We show that nonhereditary completeness for reproducing kernels
is possible in many de Branges spaces. Namely, we construct
such examples under some mild restrictions on the spectrum $\{t_n\}$
(in particular, all power growth spectra $|t_n| = |n|^\gamma$, $\gamma>0$,
$n \in \mathbb{N}$ or $n \in \mathbb{Z}$, satisfy these restrictions).

Here and later on we use the notation $U(t) \lesssim V(t)$ if there is a constant $C > 0$ such that
$U(t)\le CV(t)$ holds for all $t$ in the set in question.

\begin{theorem}[\cite{bbb}, Theorem 1.6]
\label{example}
Let $\{t_n\}$ be a sequence of real points such that
$t_n <t_{n+1}$ and $|t_n| \to \infty$\textup, $n\to\infty$. Assume that\textup,
for some $N>0$\textup,
\begin{equation}
\label{hypot}
|t_n|^{-N} \lesssim t_{n+1}-t_n =o(|t_n|) , \qquad |n|\to\infty.
\end{equation}
Then there exists a de Branges space $\he$ such that $\phi$ is of tempered growth\textup,
$\{t_n\}$ is the zero set of the function $A_E \notin \he$
and there is an exact system of reproducing kernels $\{K_\lambda\}$
in $\he$ such that its biorthogonal system is complete\textup,  but the
original system $\{K_\lambda\}$ is nonhereditarily complete.
\end{theorem}

For a question on the size of the orthogonal complement to the systems \eqref{syst3}
see Problem~5 in Section 6 below.

The statement of Theorem~\ref{example} is in sharp contrast with the fact that there exist
de Branges spaces where {\it any exact system of reproducing kernels
is hereditarily complete}. The reason for that is, however, that
a de Branges space is (in essence)
uniquely defined by the spectrum $\{t_n\}$
and the masses $\mu_n = 1/\phi'(t_n)$
of the spectral measure, while given only a spectrum there is
a lot of freedom in prescribing $\mu_n$.
Note also that in the example below condition (\ref{hypot})
is not satisfied.


\begin{example}
\label{ex-hered}
{\rm Let $\{t_n\}_{n\in \na}$ be a sequence such that $t_n>0$
and  $\inf_n t_{n+1}/t_n>1$, and let $\mu_n \equiv t_n^\gamma$, where $0<\gamma<2$.
Denote by $A_E$ the canonical product with the zeros $t_n$
and define the entire function $B_E$ by
$$
\frac{B_E(z)}{A_E(z)} = \sum_n \mu_n\bigg(\frac{1}{t_n-z} -\frac{1}{t_n}\bigg).
$$
Then the function $E = A_E-iB_E$ is in the Hermite--Biehler class.
Making use of the results of \cite[Section 5]{bbb} it is easy to show that
any exact system of reproducing kernels in $\he$ is hereditarily complete.
Otherwise, by the results of \cite{bbb}, there would exist a function $h\in \he$, $h\ne 0$, 
orthogonal to a system
\eqref{syst3} such that for a sequence $\{a_n\} \in \ell^2$ the functions
$$
\frac{h(z)}{A_E(z)}=\sum_n\frac{\overline a_n \mu_n^{1/2}}{z-t_n}\quad \text{\ and\ }
\quad \sum_n\frac{|a_n|^2}{z-t_n}
$$
have infinitely many common zeros on the real line.
However, comparing the asymptotic of the distances from these common zeros
to $\{t_n\}$, one can show that the set of common zeros can not be infinite.  }
\end{example}

Let us also mention here a result by Vasyunin \cite{va1} (see also \cite[Chapter VIII]{nik0})
characterizing the Abel summability of the Cauchy kernels
$\{(z-\overline \lambda)^{-1}:\  B(\lambda)=0\}$ expansions
in the model spaces $K_B$.


\section{Functional model for singular rank one perturbations of unbounded self-adjoint operators}



To apply the technique of the entire functions theory,
it is more convenient to work
not with perturbations of compact self-adjoint operators, but with their
unbounded inverses. These operators can be understood as
certain {\it singular} rank one perturbations of self-adjoint operators.

Let $\{t_n\}$ be an increasing sequence of real numbers (where $n \in \na$
or $n\in \zl$) such that $|t_n| \to\infty$ as $|n| \to \infty$, and let
$\mu = \sum_n \mu_n \delta_{t_n}$. In what follows $A$ is always the operator of
multiplication by the independent variable $x$ in $L^2(\mu)$ (thus, $A$
is a self-adjoint operator with simple discrete spectrum). Moreover, we
assume that $0\notin \{t_n\}$, and
so $A^{-1}$ is a bounded operator in $L^2(\mu)$.

By a {\it singular rank one perturbation} of $A$ we mean the following
operator. Let $a, b$ be two functions such that
$$
\frac{a}{x}, \  \frac{b}{x}  \in L^2(\mu)
$$
(however, possibly, $a,b \notin L^2(\mu)$).
We write $a(t_n) = a_n$, $b(t_n) = b_n$. Let $\deab \in \BC\setminus\{0\}$ be a
constant such that
\begin{equation}
\label{1+de-a-b}
\deab \ne\int_\BR x^{-1}a(x)\bar b(x) \, d\mu(x) \qquad
\text{in the case when}\ \ \ a(x) \in L^2(\mu).
\end{equation}

We associate to any such data $(a,b,\deab)$ a linear operator
$L=L(A, a,b,\deab)$, defined as follows:
$$
\begin{aligned}
 \cDL& := \big\{
y=y_0+c\cdot A^{-1}a: \\
& \qquad \qquad c\in \BC,\, y_0\in \cDA,\, \deab c+\langle y_0,
b\rangle=0
\big\};                     \\
L y& := Ay_0, \quad y\in \cDL.
\end{aligned}
$$

Condition \eqref{1+de-a-b} is equivalent to the uniqueness of the
decomposition $y=y_0+c\cdot A^{-1}a$ in the above formula for
$\cD(L)$, and so $L$ is correctly defined. Note also that $L= L(A, a, b, \deab)$
is densely defined if and only if the singular perturbation $L(A, b, a, \overline{\deab})$
is correctly defined. In this case $L^*$ is correctly defined
and $L^* = L(A, b, a, \overline{\deab})$.

This construction can be extended to what can be called
singular rank $n$ perturbations
of an unbounded operator $A$, see \cite{bar-yak}. The following fact \cite{bar-yak}
motivates these definitions: if $A$ and $L$ are
ordinary differential operators, defined by the same
regular differential expression of order $n$ and different
abstract boundary relations ($n$ independent relations in both cases),
then $L$ is a singular perturbation of $A$ of rank less or equal to $n$.

Singular rank one perturbations defined above are essentially
(unbounded) algebraic inverses to bounded rank one perturbations
of compact self-adjoint operators.
If the triplet $(a,b,\deab)$ satisfies \eqref{1+de-a-b},
then
the bounded operator
$A^{-1}- \deab^{-1} A^{-1}a (A^{-1}b)^*$
has a trivial kernel, and
%
%
\begin{align}
\label{L-Aab-kappa}
L(A, a,b, \deab)=\big(A^{-1}-\deab^{-1} A^{-1}a (A^{-1}b)^*\big)^{-1}.
\end{align}
Here we denote by $A^{-1}a (A^{-1}b)^*$ the bounded rank one operator
$A^{-1}a (A^{-1}b)^* f =  (f, A^{-1}b) A^{-1}a$, $f\in L^2(\mu)$.
Conversely, if $A^\circ$ is a compact self-adjoint operator
with a trivial kernel and $L^\circ = A^\circ + a^\circ (b^\circ)^*$
is its rank one perturbation
such that $\Ker L^\circ = 0$, then the algebraic inverse
$(L^\circ)^{-1}$ is a singular rank one perturbation of $(A^\circ)^{-1}$.

The following functional model for singular rank one perturbations
was obtained in \cite{bar-yak}. It is closely related to the model
for rank one perturbations of singular unitary operators due to Kapustin
\cite{kap} (a more general model can be found in \cite{ryzhov}).

\begin{theorem}[\cite{bar-yak}, Theorem 0.7]
\label{rank-one-model}
Let $L=L(A,a,b, \deab)$ be a singular rank one perturbation of $A$\textup,
where $b$ is a cyclic vector for the
resolvent of $A$\textup, i.e.\textup, $b_n\ne 0$ for any $n$.
Then there exists a de Branges space $\he$ and an entire function $G$
such that $E+E^* \notin \he$
\begin{equation}
\label{main0}
G\notin \he, \qquad \frac{G(z)}{z-z_0} \in \he \quad \text{if}\quad G(z_0) =0,
\end{equation}
and $L$ is unitary equivalent to the operator $T=T(E, G)$
which acts on $\he$
by the formulas
$$
\cD(T) := \{F \in \he: \text{there exists}\ c=c(F)\in \BC:
zF-cG\in \he\},
$$
$$
TF := zF - c G, \qquad F\in \cD(T).
$$

Conversely\textup, any pair $(E, G)$
where $E$ is an Hermite--Biehler function such that
$A_E \notin \he$ and $A_E(0) \ne 0$\textup, while
the entire function $G$ satisfies \eqref{main0}\textup,
corresponds to some singular rank one perturbation
$L= L(A, a, b, \deab)$ of the multiplication operator $A$ in
some space $L^2(\mu)$ with $x^{-1}  a(x)$\textup, $x^{-1}  b(x) \in L^2(\mu)$.
\end{theorem}

In fact, in \cite{bar-yak} a much more general setting is considered
where $\mu$ is an arbitrary Borel measure singular with respect to Lebesgue measure.

The functions $E=A_E-iB_E$ and $G$ appearing in the model for
$L(A,a,b, \deab)$ are related to the data $(a,b, \deab)$
by the following formulas:
\begin{equation}
\label{form-e}
\frac{B_E(z)}{A_E(z)} = \delta+ \sum_n \bigg(\frac{1}{t_n-z}
- \frac{1}{t_n}\bigg) |b_n|^2 \mu_n,
\end{equation}
where $\delta$ is an arbitrary real constant, and
\begin{equation}
\label{form-g}
\frac{G(z)}{A_E(z)} = \deab +  \sum_n \bigg(\frac{1}{t_n-z}
- \frac{1}{t_n}\bigg) a_n\overline{b_n} \mu_n,
\end{equation}
Note, in particular, that $A_E$ vanishes exactly on the set $\{t_n\}$.
The model essentially uses the expansions with respect to the orthogonal
basis $\big\{\frac{A_E}{z-t_n}\big\}$ of normalized reproducing kernels
(or, in the case of general model spaces, the representations involving Clark measures).

If we denote by $Z_G$ the zero set of $G$, then it is clear that
the eigenfunctions of $T$ are exactly the functions $\frac{G(z)}{z-\lambda}$,
$\lambda\in Z_G$. If the operator $T^*$ (equivalently, $L^*$) is well-defined,
then the eigenfunctions of $T^*$ are the functions
$\{K_\lambda\}_{\lambda \in Z_G}$. We have the following corollary.

\begin{corollary}[\cite{bar-yak}, Theorem 2.5]
\label{exact}
Let $\he$ be a de Branges space such that
$e^{i\alpha} E - e^{-i\alpha}E^* \notin \he$ for any
$\alpha \in \mathbb{R}$\textup, and let $\{t_n\}$ be the zero set of
$A_E$\textup, $t_n \ne 0$. Put $s_n = t_n^{-1}$ and let $A^\circ$
be a compact self-adjoint operator with simple eigenvalues $\{s_n\}$
and with a trivial kernel. Let
$\{K_\lambda\}_{\lambda\in \Lambda}$ be any exact system
of reproducing kernels in $\he$. Then
there exists a bounded rank-one perturbation $L^\circ$ of $A^\circ$
such that
$$
\sigma(L^\circ)=\big\{\la^{-1}: \la\in \Lambda\big\},
$$
$\ker L^\circ=\ker (L^\circ)^*=0$,
and there is a unitary transform $U:\he\to H$\textup,
which maps the system $\{K_\lambda\}_{\lambda\in \Lambda}$
to a system of eigenvectors of $L^\circ$\textup:
$L^\circ\big(U K_\lambda\big) = \la^{-1}\big(U K_\lambda\big)$.
\end{corollary}

It follows that at the same time, $U^*$ takes the
biorthogonal system to $\{K_\lambda\}_{\lambda\in \Lambda}$  into a
system of eigenvectors of the adjoint operator $(L^\circ)^*$.



\section{Applications of the functional model}

Using the above functional model we obtain a number of results
concerning the following natural questions:
\smallskip

($i$) {\it When does $L = L(A, a, b, \deab)$ have a complete system
of eigenvectors \textup(i.e.\textup, $L$ is complete\textup)\textup?}
\smallskip

($ii$) {\it When does the completeness of $L$ imply the completeness of $L^*$\textup?}
\smallskip

$(iii)$  {\it When does $L$ admit the spectral synthesis\textup?}
\smallskip

$(iv)$ {\it For which $A$ does there exist a rank one perturbation
$L(A, a, b, \deab)$ with the spectrum at infinity
\textup(i.e.\textup, $L$ is the inverse to a Volterra operator\textup)\textup? }

\subsection{Completeness of rank one perturbations of self-adjoint operators}
We say that $L(A, a, b, \deab)$ is a {\it generalized weak
perturbation}, if
\begin{equation}
\label{Cond-pos}
\sum_n \frac{|a_n b_n| \mu_n}{|t_n|}<\infty, \qquad
\sum_n \frac{a_n\overline b_n \mu_n}{t_n} \ne \deab.
\end{equation}
It is essentially a corollary of Matsaev's theorem
\cite{Mats61} that any generalized weak perturbation is complete.
We get another sufficient condition for completeness of $L$
if we assume a certain positivity condition.

\begin{theorem}[\cite{bar-yak}, Theorems 0.1, 3.3]
\label{positive}
1. If $L=L(A,a,b,\deab)$ is a generalized weak perturbation\textup, then
$L^*$ is correctly defined\textup, and $L$ and $L^*$ are complete.

2. Suppose that $a_n\overline b_n\ge 0$ for all but
a finite number of values of $n$ and $\sum_n |t_n|^{-1}|a_nb_n|\mu_n = \infty$.
Then $L^*$ is correctly defined\textup, and $L$ and $L^*$ are complete.
\end{theorem}

A typical example when (\ref{Cond-pos}) is satisfied is
that there is $\alpha\in [0,1]$ such that
$$
a \in |x|^\alpha L^2(\mu), \qquad b\in |x|^{1-\alpha} L^2(\mu).
$$
At the same time it is easy to show that even for rank one perturbations,
when we relax slightly the generalized weakness property \eqref{Cond-pos},
the resulting perturbation may become the inverse to a Volterra operator.
We state the corresponding result for perturbations of compact operators.

\begin{theorem}[\cite{bar-yak}, Theorem 0.6]
\label{sharp}
There exists a sequence $s_n \to 0$ and a measure $\mu = \sum_n \mu_n
\delta_{s_n}$
with the following property\textup: for any $\alpha_1, \alpha_2 \ge 0$ with
$\alpha_1+ \alpha_2 <1 $
there exist $a^\circ \in |x|^{\alpha_1} L^2(\mu)$ and $b^\circ \in
|x|^{\alpha_2} L^2(\mu)$
such that the perturbed operator $A^\circ + a^\circ (b^\circ)^*$  \textup(where
$A^\circ$  is the operator
of multiplication by $x$ in $L^2(\mu)$\textup)  is a Volterra operator.
\end{theorem}

Positive results on the completeness can also be obtained
for rank one perturbations of $A$ that are not
generalized weak perturbations, under the assumption that the spectrum
of $A$ is \textit{exponentially sparse}, which means that
$s_{n+1}/s_n<\gamma$ for all $n>0$ (with a similar condition for
$n<0$) and for some small (absolute) constant $\gamma\in(0,1)$.


\subsection{Relations between completeness
of the perturbed operator and of its adjoint}
If a bounded operator $T$ is complete, a trivial obstacle
for completeness of $T^*$ is that $T$ may have a nontrivial kernel, while
$\ker T^* = 0$. The first (highly nontrivial) examples
of the situation where a compact operator
$T$ is complete and $\ker T=0$, while $T^*$
is not complete, were constructed by Hamburger \cite{hamb}.
In \cite{DeckFoPea} Deckard, Foia\c{s}, and Pearcy gave a simpler
construction.
By definition, in these examples, the eigenvector system of $T$ cannot be hereditarily
complete, because if it were true, the biorthogonal to this system (which are just
the eigenvectors of $T^*$) would be also complete.

In the above-mentioned examples one cannot conclude that the corresponding operator is a small
(for instance, finite rank or trace class)
perturbation of a self-adjoint operator.
Corollary \ref{adjoint} below shows that one can find such examples
among rank one perturbations of compact self-adjoint operators
with almost arbitrary spectrum.

We start with the situation of singular rank one perturbations
of unbounded operators. As follows from the functional model, completeness
of $L$ and $L^*$ is reduced to the completeness of a system
of reproducing kernels and of its biorthogonal system.
The relations between these two completeness problems were studied in \cite{bb}.

The next theorem shows that under certain additional assumptions
the completeness of a singular rank one perturbation
$L$ implies the  completeness of its adjoint.

\begin{theorem}[\cite{bar-yak}, Theorem 0.2]
\label{frombb}
Let the data $(a, b, \deab)$ satisfy the property $a\notin L^2(\mu)$\textup,
and let the perturbation $L=L(A, a,b, \deab)$ be complete.
Assume that its adjoint $L^*$ is  correctly defined. Then $L^*$
is also complete if any of the following conditions is fulfilled\textup:

$(i)$ $|a_n|^2 \mu_n  \lesssim |t_n|^{-N}$ for some $N>0$\textup;

$(ii)$ $|b_n a_n^{-1}| \lesssim |t_n|^{-N}$ for some $N>0$.
\end{theorem}

In general, the situation is much subtler.
Applying the results and methods from \cite{bb} we
are able to give examples
when the adjoint to a complete perturbation fails
to be complete.

\begin{theorem}[\cite{bar-yak}, Theorem 0.3]
\label{noncompleteness2}
For any cyclic self-adjoint operator $A$ with
discrete spectrum\textup, 
there exists a singular rank one
perturbation $L$ of $A$ with real spectrum\textup,
which is not complete\textup, while its adjoint $L^*$ is correctly defined\textup,
has trivial kernel and is complete.
Moreover\textup, the orthogonal complement to the space spanned by the
eigenvectors of $L$ may be infinite-dimensional.
\end{theorem}

\begin{corollary}[\cite{bar-yak}, Corollary 0.4]
\label{adjoint}
For any compact self-adjoint operator $A^\circ$ with
simple spectrum $\{s_n\}$\textup, $s_n \ne 0$\textup,
there exists a bounded rank one perturbation $L^\circ$ of $A^\circ$
with real spectrum such that $L^\circ$ is complete
and $\ker L^\circ = 0$\textup, while $(L^\circ)^*$ is not complete.
\end{corollary}


\subsection{Spectral synthesis for rank one perturbations}

In view of the theorem of Markus \cite[Theorem 4.1]{markus}
and Corollary \ref{exact}, the results of Section 3 may be interpreted
as the results about the spectral synthesis for rank one perturbations
of compact self-adjoint operators. In particular,
we give now the accurate statement for Theorem \ref{synthesis-1}.

\begin{theorem}[\cite{bar-yak}, Theorem 0.5]
\label{synthesis}
Let $\{s_n\}$ be a sequence of real numbers
\textup(ordered so that $s_n>0$ and $s_n$ decrease for
$n\ge 0$\textup, and $s_n<0$ and increase for $n<0$\textup) and
assume that for some $N>0$
$$
|s_n|^N \lesssim |s_{n+1} - s_n| = o(|s_n|), \qquad |n| \to \infty.
$$
Let $A^\circ$ be a compact self-adjoint operator with simple
eigenvalues $\{s_n\}$ and with a trivial kernel.
Then there exists a rank one perturbation $L^\circ$ of $A^\circ$
with real spectrum such that both $L^\circ$
and $(L^\circ)^*$ have complete sets of eigenvectors\textup, but
$L^\circ$ does not admit the spectral synthesis.
\end{theorem}

As Marcus showed in \cite{markus},
most classical sufficient conditions of completeness
imply already that the operator admits the spectral synthesis.


\subsection{Removability of the spectrum}

In this subsection we address the following question, which, in a sense,
is opposite to the completeness problem:
\textit{For which measures $\mu = \sum_n \mu_n \delta_{t_n}$
does there exist a singular perturbation $L$ of $A$\textup, whose spectrum is empty
\textup(in other words, consists only of the point at infinity\textup)\textup?}

Thus, the problem is to describe those spectra $\{t_n\}$
for which the spectrum of the perturbation is empty.
Such spectra will be said to be {\it removable}. It is clear that the
property to be removable or nonremovable  depends only on $\{t_n\}$, but
not on the choice of the masses $\mu_n$.

The removability criterion will be given
in terms of entire functions of the
so-called {\it Krein class}. We say that an entire function
$F$ is in the Krein class $\mathcal{K}_1$, if it has only real simple zeros
$t_n$ and can be represented as
$$
\frac{1}{F(z)} = q+ \sum_n c_n
\Big(\frac{1}{t_n-z}-\frac{1}{t_n}\Big), \qquad \sum _n
t_n^{-2}|c_n| <\infty,
$$
where $c_n = -1/F'(t_n)$ and $q=1/F(0)$.

Then our main result in this section reads as follows:

\begin{theorem}[\cite{bar-yak1}, Theorem 0.1]
\label{annih2}
Let $t_n \in \RR$ and $|t_n| \to\infty$\textup, $|n| \to \infty$.
The following are equivalent\textup:

$(i)$ The spectrum $\{t_n\}$ is removable\textup;

$(ii)$ There exists a function $F \in \mathcal{K}_1$
%
%
whose zero set coincides with $\{t_n\}$.
\end{theorem}

A somewhat unexpected and counterintuitive
consequence of Theorem \ref{annih2} is that
adding a finite number of points to the spectrum helps it to
become removable, while deleting  a finite number of points can
make it nonremovable. E.g., the spectrum $t_n = n+1/2$,
$n\in\mathbb{Z}$, is removable, while $t_n = n+1/2$, $n\in \mathbb{Z}$,
$n\ne 0$, is not. Also, the spectrum
$\{n^2\}_{n\in \mathbb{N}}$ is removable, but
$\{n^2\}_{n\ge 2}$ is nonremovable.

In view of the relation between singular rank one perturbations
(see \eqref{L-Aab-kappa})
and usual rank one perturbations of bounded self-adjoint operators,
we have an immediate counterpart of Theorem  \ref{annih2}
for Volterra rank one perturbation of compact operators.

\begin{theorem}[\cite{bar-yak1}, Theorem 0.2]
\label{annih3}
Let $s_n \in \RR$\textup, $s_n \ne 0$\textup,
and $|s_n| \to 0$\textup, $|n| \to \infty$\textup, and let
$A^\circ$ be a compact self-adjoint operator with simple
point spectrum $\{s_n\}$ and with a trivial kernel.
The following are equivalent\textup:

$(i)$ There exists a rank one perturbation $L^\circ$ of
$A^\circ$ which is a Volterra operator\textup;

$(ii)$ The points $t_n = s_n^{-1}$ form the zero set of some function
$F \in \mathcal{K}_1$.
\end{theorem}

We find it a bit surprising that such a natural question was not
previously addressed. There is a vast literature on the subject, and many results
(mostly due to Krein, Gohberg, and Matsaev) concern the relations between
the real and the imaginary parts of a Volterra operator;
see \cite[Chapter III]{gkvol} and
\cite[Chapter IV]{Gohb_Krein} (especially, Section 10, where some
partial results about Volterra operators with finite rank imaginary
parts are obtained). Let us also mention a beautiful theorem of Livshits
which says that a dissipative rank one perturbation of a self-adjoint operator
is unitary equivalent to the integration operator
(see \cite{liv} or \cite[Chapter 1, Theorem 8.1]{gkvol}).
Still, we were unable to find any result
explicitly describing the compact self-adjoint operators
whose rank one perturbation is a Volterra operator.

To conclude, we deduce the above-mentioned theorem of Livshits
from the functional model of rank one
perturbations\footnote{It was N. Nikolski who attracted
our attention to the Livshits theorem and suggested to deduce it
using our methods.}.
Namely, we show the following:

\begin{theorem}[Livshits, \cite{liv}]
\label{livs}
Let $L^\circ=A^\circ +iB^\circ$
be a Volterra operator \textup(in some Hilbert space $H$\textup)
where both $A^\circ$ and $B^\circ$
are self-adjoint and $B^\circ$ is of rank one.
Then the point spectrum of $A^\circ$ is given by $s_n = c (n+1/2)^{-1}$\textup,
$n\in \mathbb{Z}$\textup, for some $c>0$.
\end{theorem}

From this, one can deduce that $A^\circ$ is unitary
equivalent to the integral operator (having the same spectrum)
$$
(\tilde A f) (x) = i\int_0^{2\pi c} f(t)\, \sig (x-t)\,dt, \qquad
f\in L^2(0,2\pi c),
$$
while $L^\circ$ is unitary equivalent to the integration operator
$(\tilde L f) (x) = 2i \int_0^{2\pi c} f(t)\,dt$.

Since $B^\circ=(B^\circ)^*$, we have $B^\circ x =  (x, b^\circ) b^\circ$
for some $b^\circ \in H$.
Passing to the unbounded inverses we obtain
(after an obvious unitary equivalence) a singular rank one perturbation
$L = L(A, a, b, \deab)$
of the operator $A$ of multiplication by the independent variable
in some space $L^2(\mu)$ where $\mu = \sum_n \mu_n \delta_{t_n}$, $t_n
= s_n^{-1}$. Moreover, in the case of the self-adjoint imaginary part,
we may assume that $\deab = -1$ and  $a = i b$.

Applying the functional
model from Section 4 we construct a de Branges space $\he$ and a function $G$
as in (\ref{form-e})--(\ref{form-g}). In our case we have
$$
\begin{aligned}
\frac{B_E(z)}{A_E(z)} & = \delta+ \sum_n \bigg(\frac{1}{t_n-z}
- \frac{1}{t_n}\bigg) |b_n|^2 \mu_n, \\
\frac{G(z)}{A_E(z)} & = -1 +  i \sum_n \bigg(\frac{1}{t_n-z}
- \frac{1}{t_n}\bigg) |b_n|^2 \mu_n,
\end{aligned}
$$
whence $G = -A_E + i(B_E -\delta A_E)$.

Since $L$ (and, thus, the model operator $T$)
is the inverse to a Volterra operator, the spectrum of $T$ is
the point at infinity. Thus, $G$ has no zeros in $\mathbb{C}$. Also, by the
results of \cite{bar-yak1}, the function $E$ is of Cartwright class
and the same is true for $G$. We conclude that $G(z) = \exp(i \pi cz)$
for some real $c$. Thus,
$$
e^{i\pi cz} =  -A_E(z) + i\big(B_E(z) -\delta A_E(z)\big).
$$
The functions $A_E$ and $B_E$ are real on the real axis. Taking the real parts,
we have $A_E(z) = -\cos \pi cz$, and so $t_n = c^{-1} (n+1/2)$,
$n\in \mathbb{Z}$, as required.


\section{Open problems}

Here we mention a few questions related to the above results.
A very basic question, which remains open, is the following one.
\medskip
\\
{\bf Problem 1.} { \it Is it true that any hereditarily complete system
of exponentials in $L^2(-\pi,\pi)$ is a linear summation basis \textup(i.e.\textup,
the corresponding Fourier series
are all summable by a linear summation method\textup)\textup?}
\medskip

We are able to construct systems $\{K_\lambda\}_{\lambda\in \Lambda}$
of reproducing kernels in a Paley--Wiener space such that
for some partition $\Lambda = \Lambda_1 \cup \Lambda_2$
the orthogonal complement in $\pw$ to the system
$\mathcal F_{\Lambda_1,\Lambda_2}=\{g_\lambda\}_{\lambda\in \Lambda_1} \cup \{K_\lambda\}_{\lambda\in \Lambda_2}$
is one-dimensional.
\medskip
\\
{\bf Problem 2.} {\it How to characterize the vectors orthogonal to
such systems\textup?
}
\medskip

It is clear that given a partition
$\Lambda = \Lambda_1 \cup \Lambda_2$ such that $\mathcal F_{\Lambda_1,\Lambda_2}$ has
a nontrivial orthogonal complement, the orthogonal
complement to the system $\mathcal F_{\Lambda_1 \setminus \{\lambda_0\},\Lambda_2 \cup\{\lambda_0\}}$
also will be nontrivial. Thus, we can
always move a finite number of points from $\Lambda_1$ to $\Lambda_2$ and in the opposite direction.
\medskip
\\
{\bf Problem 3.} {\it What is the structure of the set of those partitions
$\Lambda = \Lambda_1 \cup \Lambda_2$ for which the orthogonal complement
to $\mathcal F_{\Lambda_1,\Lambda_2}$
is nontrivial\textup? Is this set in a sense \textup"connected\textup"\textup?}
\medskip

As we have seen in Example \ref{ex-hered}, there exist de Branges spaces
$\he$ where any exact system of reproducing kernels is hereditarily complete.
The proof in this example uses essentially that the spectrum $\{t_n\}$
is lacunary. However, a de Branges space is uniquely (up to a natural isometry) defined by
its spectral data $t_n$ and $\mu_n = 1/\phi'(t_n)$. Apart from the case of
the Paley--Wiener spaces, in our examples of nonhereditarily complete systems only $t_n$ were fixed.
\medskip
\\
{\bf Problem 4.} {\it Describe spectral data $(t_n, \mu_n)$ such that any exact system
of the reproducing kernels in the corresponding de Branges space $\he$ is
hereditarily complete.  Is it true that for the spectrum $\{t_n\}$ satisfying
\eqref{hypot} and for any $\mu_n$ with $\sum_n \mu_n (t_n^2+1)^{-1}<\infty$
there exists a nonhereditarily complete system in $\he$\textup? }
\medskip

It should be mentioned that all our examples of nonhereditary completeness
are constructed in the {\it reverse order}. Namely, we start with a vector
$h$ (in $\pw$ or $\he$) with some special
properties and then construct a sequence $\Lambda$
such that $h$ is orthogonal to some system $\mathcal F_{\Lambda_1,\Lambda_2}$.
Therefore, we are not able to produce a nonhereditarily complete system
of reproducing kernels in a de Branges space such that
the orthogonal complement to the system $\mathcal F_{\Lambda_1,\Lambda_2}$
is two-dimensional
or infinite-dimensional. The existence of such examples is an intriguing problem.
Note that it follows from (the proof of) the theorem of Markus \cite[Theorem 4.1]{markus}
that if the system of eigenvectors is hereditarily complete
up to a finite-dimensional defect, then the spectral synthesis holds up
to a defect of the same
\medskip
dimension.
\\
{\bf Problem 5.} {\it Do there exist nonhereditarily complete
systems of reproducing kernels in a de Branges space such that for some
partition the orthogonal complement to $\mathcal F_{\Lambda_1,\Lambda_2}$
is $n$-dimensional
with $n\ge 2$ or even infinite-dimensional\textup?\,\footnote{
{\bf Added in proof.} Recently the first three authors have constructed an example
of a de Branges space with a nonhereditarily complete system
$\{K_\lambda\}_{\lambda\in \Lambda}$ such that, for some partition,
the system $\mathcal F_{\Lambda_1,\Lambda_2}$ has infinite defect.}}
\medskip

Finally, though the spectral synthesis may fail even for exponential systems
(reproducing kernels in the Paley--Wiener space), it may of interest  to
distinguish those systems for which it still holds. The generating function $G$
may provide a natural language for this problem.
\medskip
\\
{\bf Problem 6.} {\it What conditions on the generating function $G$ of an exact
system of reproducing kernels in $\he$ ensure that the system is hereditarily complete\textup?}
\medskip

Recently, this problem was considered by Gubreev and Tarasenko \cite{gub}
who showed that in the case when $|G/E|^2$ is a Muchenhoupt $A_2$-weight
on $\mathbb{R}$, the corresponding system of kernels is hereditarily complete.
\medskip

To conclude, it seems that hereditary completeness for systems of reproducing kernels
in de Branges spaces and the related spectral synthesis problems
(even for such simple class of linear operators
as rank one perturbations of self-adjoint operators)
remain a rich field where there is still much to explain and to explore.

\end{document}